\documentclass[fleqn, a4paper]{article}

\usepackage{amsmath,txfonts,bm}
\usepackage{graphicx}
\usepackage{theorem}
\theorembodyfont{\upshape}
\theoremstyle{break}
\newtheorem{thm}{Theorem}[section]
\newtheorem{lem}[thm]{Lemma}
\newtheorem{rmk}[thm]{Remark}
\newtheorem{cor}[thm]{Corollary}
\usepackage{algorithm,algorithmic}
\usepackage{ascmac}

\begin{document}

    \noindent 
    \LARGE{A new subtraction-free formula for lower bounds of the minimal singular value of an upper bidiagonal matrix}  \\

    \large{Takumi Yamashita
          \footnote{28-20 Kojogaoka, Otsu, Shiga 520-0821 Japan, Tel.: +81-77-522-7447, Fax:  +81-77-522-7447, e-mail: t-yamashita@kke.biglobe.ne.jp},
          Kinji Kimura
          \footnote{Graduate School of Informatics, Kyoto University, Yoshida-Honmachi, Sakyo-ku, Kyoto, Kyoto 606-8501 Japan, e-mail: kkimur@amp.i.kyoto-u.ac.jp}
          and
          Yusaku Yamamoto
          \footnote{Graduate School of Informatics and Engineering, The University of Electro-Communications, 1-5-1, Chofugaoka, Chofu, Tokyo 182-8585, Japan,
                    e-mail: yusaku.yamamoto@uec.ac.jp}
          }

    \section*{Abstract}

    Traces of inverse powers of a positive definite symmetric tridiagonal matrix give lower bounds of the minimal singular value of an upper bidiagonal matrix.
    In a preceding work, a formula for the traces which gives the diagonal entries of the inverse powers is presented.
    In this paper, we present another formula which gives the traces based on a quite different idea from the one in the preceding work.
    An efficient implementation of the formula for practice is also presented.

    \section{Introduction}  \label{Intro}

    A lower bound of the minimal singular value of a matrix has historically been investigated for estimation of an upper bound of the condition number of a matrix.
    As another application, such a lower bound for an upper bidiagonal matrix may be used to accelerate convergence of iteration in some singular value computing algorithms
    \cite{FP94, IN06, vonMatt97, NIN12}.
    In the standard procedure for computing the singular values,
    one first reduces the input matrix to an upper bidiagonal matrix by orthogonal transformations
    and then computes the singular values of the obtained upper bidiagonal matrix by some iterative algorithms.
    The iterative algorithms referred above use a technique called the shift of origin.
    This technique requires a quantity called a shift.
    A lower bound of the minimal singular value of the upper bidiagonal matrix can be used to determine this quantity.

    Several lower bounds of the minimal singular value of a matrix have been proposed.
    For example, see \cite{HP92, Johnson89, JS98, Li01, Rojo99, YG97}.
    For an upper bidiagonal matrix $B$, where all the diagonal and the upper subdiagonal entries are positive,
    the traces $\textrm{Tr} ( ( B B^{\top} )^{- M} )$ $( M = 1, 2, \dots )$ give lower bounds of the minimal singular value of $B$.
    For example, the following two quantities
    \begin{equation*}
        \varrho = ( \textrm{Tr} ( ( B B^{\top} )^{- 1} ) )^{- \frac{1}{2}}
    \end{equation*}
    and
    \begin{equation*}
        \upsilon = \sqrt{\frac{1}{\textrm{Tr} ( ( B B^{\top} )^{- 1} )}}
                   \cdot \sqrt{\frac{N}{1 + \sqrt{( N - 1 ) \left( N \cdot \displaystyle \frac{  \textrm{Tr} ( ( B B^{\top} )^{- 2} )      }
                                                                                              {( \textrm{Tr} ( ( B B^{\top} )^{- 1} ) )^{2}} - 1 \right)}}},
    \end{equation*}
    where $N$ is the matrix size of $B$, are such lower bounds.
    For details, see \cite{vonMatt97} by von Matt.
    For computation of the traces of $( B B^{\top} )^{- 1}$ and $( B B^{\top} )^{- 2}$,
    von Matt \cite{vonMatt97} also presented a method to compute the diagonal entries of these inverses.
    On the other hand, Kimura {\it et al}. \cite{KYN11} presented a sequence of lower bounds of the minimal singular value of $B$.
    These lower bounds $\theta _{M} ( B )$ $( M = 1, 2, \dots )$ are given with the traces $J_{M} ( B ) = \textrm{Tr} ( ( B^{\top} B )^{- M} ) = \textrm{Tr} ( ( B B^{\top} )^{- M} )$ as
    \begin{equation*}
        \theta _{M} ( B ) = ( J_{M} ( B ) )^{- \frac{1}{2 M}},  \qquad  M = 1, 2, \dots .
    \end{equation*}
    It holds $\varrho = \theta _{1} ( B )$.
    They increase monotonically and converge to the minimal singular value $\sigma _{\min} ( B )$ of $B$ as $M$ goes to infinity \cite[Theorem 3.1]{KYN11}, that is,
    \begin{gather*}
        \theta _{1} ( B ) < \theta _{2} ( B ) < \cdots < \sigma _{\min} ( B ),  \\
        \lim_{M \to \infty} \theta _{M} ( B ) = \sigma _{\min} ( B ).
    \end{gather*}
    Kimura {\it et al}. \cite{KYN11} also presented a formula for computation of the traces of $J_{M} ( B )$ for an arbitrary positive integer $M$.
    This formula gives the diagonal entries of the inverse powers $( B^{\top} B )^{- M}$ and $( B B^{\top} )^{- M}$ $( M = 1, 2, \dots )$ in a form of recurrence relation.
    In \cite{YKN12}, Yamashita {\it et al}. derived another formula for these diagonal entries starting from the formula in \cite{KYN11}.
    While the formula in \cite{KYN11} includes subtraction in it in the case of $M \geq 2$,
    the formula in \cite{YKN12} consists of only addition, multiplication and division among positive quantities.
    Namely, the formula in \cite{YKN12} is ``subtraction-free''.
    This property clearly excludes any possibility of cancellation error.
    
    In this paper, we present another formula for computation of the traces $J_{M} ( B )$ $( M = 1, 2, \dots )$.
    This formula is also subtraction-free.
    We derive the formula with an idea which is quite different from that in \cite{YKN12}.
    We do not aim to obtain the diagonal entries of $( B^{\top} B )^{- M}$ or $( B B^{\top} )^{- M}$ $( M = 1, 2, \dots )$ in the derivation.
    Instead, equations on the determinant and the entries of $A - \lambda I$, where $A$ is    $B^{\top} B$
                                                                                           or $B B^{\top}$, $\lambda$ is a parameter and $I$ is the unit matrix, are considered.
    The new formula is obtained by differentiating these equations with respect to the parameter $\lambda$ repeatedly.
    Computational cost for the traces are also discussed.
    Moreover, an implementation for the trace $J_{2} ( B )$ which is useful in practice is presented.
    This implementation has the following merits compared with that in \cite{YKN12}.
    \begin{itemize}
        \item The number of operations is smaller compared with the implementation in \cite{YKN12}.
        \item Only one ``loop'' is required while the implementation in \cite{YKN12} requires two loops.
        \item No ``array'' is necessary, in contrast to the implementation in \cite{YKN12}.
    \end{itemize}

    This paper is organized as follows.
    In Section \ref{Derivation}, the new formula is derived.
    In Section \ref{Cost}, computational cost for the new formula is discussed.
    In Section \ref{Implementation}, an efficient implementation for the trace $J_{2} ( B )$ is presented.
    Section \ref{ConcRmks} is devoted for concluding remarks.

    \section{Derivation of the formula for the traces}  \label{Derivation}

    Let us consider an $N \times N$ real upper bidiagonal matrix $B$, where all the diagonal and the upper subdiagonal entries are positive.
    In this section, we derive the new formula for the traces $J_{p} ( B ) = \textrm{Tr} ( ( B^{\top} B )^{- p} ) = \textrm{Tr} ( ( B B^{\top} )^{- p} )$
                     for an arbitrary positive integer $p$ in a form of recurrence relation.
    From these traces, lower bounds of the minimal singular value of $B$ are obtained.
    In the context of singular value computation, we can assume the positivity of the diagonal and the upper subdiagonal entries of $B$ without loss of generality \cite{FP94}.
    We present two recurrence relations in Sections \ref{Sect_DerivI} and \ref{Sect_DerivII}.
    The ideas to derivate the recurrence relations are quite different from those in \cite{YKN12}.

    Hereafter, we fix some notations.
    Let $B$ be
    \begin{equation}  \label{UpperBidiagonal}
        B = \left( \begin{array}{ccccc}
                       \sqrt{q_{1}}  &  \sqrt{e_{1}}  &                 &                    &                    \\
                                     &  \sqrt{q_{2}}  &  \sqrt{e_{2}}   &                    &                    \\
                                     &                &  \ddots         &  \ddots            &                    \\
                                     &                &                 &  \sqrt{q_{N - 1}}  &  \sqrt{e_{N - 1}}  \\
                                     &                &                 &                    &  \sqrt{q_{N    }}
                   \end{array}                                                                                       \right) ,
    \end{equation}
    where $q_{i} > 0$ for $i = 1, \dots , N$ and $e_{i} > 0$ for $i = 1, \dots , N - 1$.
    Let $I$ be the $N \times N$ unit matrix.
    We use the convention $\sum_{i = j}^{k} = 0$ if $j > k$. 
    Let $\lambda$ be a parameter.

    \subsection{Derivation - type I}  \label{Sect_DerivI}
    
    In this subsection, we actually derive the formula.
    Let the eigenvalues of $B^{\top} B$ be $\lambda _{1}, \dots , \lambda _{N}$.
    For an arbitrary positive integer $p$, the eigenvalues of $( B^{\top} B )^{- p}$ are $\lambda _{1}^{- p}, \dots , \lambda _{N}^{- p}$.
    Then, the summation $\sum_{i = 1}^{N} \lambda _{i}^{- p}$ is  the trace of $( B^{\top} B )^{- p}$.
    We derive a formula to compute this summation.
    The matrix $B^{\top} B$ is given as
    \begin{equation*}
        B^{\top} B = \left( \begin{array}{cccc}
                                      q_{1}         &  \sqrt{q_{1}   e_{1}}  &                              &                                \\
                                \sqrt{q_{1} e_{1}}  &        q_{2} + e_{1}   &  \ddots                      &                                \\
                                                    &  \ddots                &  \ddots                      &  \sqrt{q_{N - 1}   e_{N - 1}}  \\
                                                    &                        &  \sqrt{q_{N - 1} e_{N - 1}}  &        q_{N    } + e_{N - 1}
                            \end{array}                                                                                                          \right) .
    \end{equation*}
    It can be readily verified that we obtain the following matrix
    \begin{equation*}
        A = \left( \begin{array}{cccc}
                       q_{1}  &  q_{1}   e_{1}  &          &                         \\
                       1      &  q_{2} + e_{1}  &  \ddots  &                         \\
                              &  \ddots         &  \ddots  &  q_{N - 1}   e_{N - 1}  \\
                              &                 &  1       &  q_{N    } + e_{N - 1}
                   \end{array}                                                           \right)
    \end{equation*}
    by similarity transformation.
    Then, the matrices $A$ and $B^{\top} B$ have the same eigenvalues.
    A key point of this derivation is to express the determinant of $A - \lambda I$ in two ways.
    As the first way, the determinant is expressed as
    \begin{equation}  \label{DetAlIEigenv}
        \det ( A - \lambda I ) = \prod_{i = 1}^{N} ( \lambda _{i} - \lambda ).
    \end{equation}
    For the second way, let us consider decomposition of the matrix
    \begin{equation*}
        A - \lambda I = \left( \begin{array}{cccc}
                                   q_{1} - \lambda  &  q_{1}   e_{1}            &          &                                   \\
                                   1                &  q_{2} + e_{1} - \lambda  &  \ddots  &                                   \\
                                                    &  \ddots                   &  \ddots  &  q_{N - 1}   e_{N - 1}            \\
                                                    &                           &  1       &  q_{N    } + e_{N - 1} - \lambda
                               \end{array}                                                                                        \right)
    \end{equation*}
    into the matrix product expressed as
    \begin{equation*}
        A - \lambda I = \left( \begin{array}{cccc}
                                    \hat{q}_{1}^{( 0 )}  &                       &          &                       \\
                                    1                    &  \hat{q}_{2}^{( 0 )}  &          &                       \\
                                                         &  \ddots               &  \ddots  &                       \\
                                                         &                       &  1       &  \hat{q}_{N}^{( 0 )}
                                \end{array}                                                                            \right)
                         \left( \begin{array}{cccc}
                                    1  &  \hat{e}_{1}^{( 0 )}  &          &                           \\
                                       &  1                    &  \ddots  &                           \\
                                       &                       &  \ddots  &  \hat{e}_{N - 1}^{( 0 )}  \\
                                       &                       &          &  1
                                \end{array}                                                                            \right) ,
    \end{equation*}
    where $\hat{q}_{i}^{( 0 )}$ for $i = 1, \dots , N$ and $\hat{e}_{i}^{( 0 )}$ for $i = 1, \dots , N - 1$ are functions of $\lambda$.
    These functions are repeatedly differentiated in the discussion shown below.
    The superscript $( 0 )$ indicates that the function has not been differentiated yet.
    Comparison of the diagonal and the upper subdiagonal entries of $A - \lambda I$ gives
    \begin{align}
        &q_{i} + e_{i - 1} - \lambda = \hat{q}_{i}^{( 0 )} + \hat{e}_{i - 1}^{( 0 )},  & &i = 1, \dots , N    ,  \label{eq_hat_sum}   \\
        &q_{i}   e_{i    }           = \hat{q}_{i}^{( 0 )}   \hat{e}_{i    }^{( 0 )},  & &i = 1, \dots , N - 1,  \label{eq_hat_prod}
    \end{align}
    where $e_{0}^{(0)} = 0$ and $\hat{e}_{0}^{(0)} = 0$.
    Then, the functions $\hat{q}_{i}^{( 0 )}$ for $i = 1, \dots , N$ and $\hat{e}_{i}^{( 0 )}$ for $i = 1, \dots , N - 1$ are obtained by the following recurrence relation
    \begin{align*}
        &\hat{q}_{1}^{( 0 )} =                     q_{1}             - \lambda                           ,  & &                         \\
        &\hat{e}_{i}^{( 0 )} = \displaystyle \frac{q_{i}   e_{i    }}{           \hat{q}_{i    }^{( 0 )}},  & &  i = 1, \dots , N - 1,  \\
        &\hat{q}_{i}^{( 0 )} =                     q_{i} + e_{i - 1} - \lambda - \hat{e}_{i - 1}^{( 0 )} ,  & &  i = 2, \dots , N.
    \end{align*}
    Thus, the second expression of the determinant of $A - \lambda I$ is given as
    \begin{equation}  \label{DetAlIq0}
        \det ( A - \lambda I ) = \prod_{i = 1}^{N} \hat{q}_{i}^{( 0 )}.
    \end{equation}
    By (\ref{DetAlIEigenv}) and (\ref{DetAlIq0}), we have
    \begin{equation}  \label{EqDetAlI}
        \prod_{i = 1}^{N} ( \lambda _{i} - \lambda ) = \prod_{i = 1}^{N} \hat{q}_{i}^{( 0 )}.
    \end{equation}
    We differentiate this equation (\ref{EqDetAlI}).
    The result of differentiation of the left-hand side of (\ref{EqDetAlI}) is
    \begin{equation}  \label{diff_lambda}
                   \sum_{i = 1}^{N} \left( - \frac{1}{\lambda _{i} - \lambda} \prod_{j = 1}^{N} ( \lambda _{j} - \lambda ) \right)
        = \left( - \sum_{i = 1}^{N}          \frac{1}{\lambda _{i} - \lambda}                                                      \right) \det ( A - \lambda I )  .
    \end{equation}
    Before differentiation of the right-hand-side of (\ref{EqDetAlI}), we introduce functions $\hat{q}_{i}^{( p )}$ of $\lambda$ defined by
    \begin{equation}  \label{def_hat_q_p}
        \hat{q}_{i}^{( p )} = \frac{d^{p} \hat{q}_{i}^{( 0 )}}{d \lambda ^{p}}
    \end{equation}
    for $i = 1, \dots , N$ and $p = 1, 2, \dots$.
    The result of differentiation is
    \begin{equation}  \label{diff_q}
                 \sum_{i = 1}^{N} \frac{\prod_{j = 1}^{N} \hat{q}_{j}^{( 0 )}}{\hat{q}_{i}^{( 0 )}} \cdot \hat{q}_{i}^{( 1 )}
        = \left( \sum_{i = 1}^{N} \frac{                  \hat{q}_{i}^{( 1 )}}{\hat{q}_{i}^{( 0 )}}                           \right) \det ( A - \lambda I ).
    \end{equation}
    From (\ref{diff_lambda}) and (\ref{diff_q}), we derive
    \begin{equation}  \label{Eq_diff_1}
          \sum_{i = 1}^{N}          \frac{1                  }{     \lambda _{i}         - \lambda}
        = \sum_{i = 1}^{N} \left( - \frac{\hat{q}_{i}^{( 1 )}}{\hat{q      }_{i}^{( 0 )}          } \right) .
    \end{equation}
    Then, by substituting $\lambda = 0$ into the left-hand-side of (\ref{Eq_diff_1}),
          we have the summation $\sum_{i = 1}^{N} \lambda _{i}^{- 1}$ which is equal to the trace of $( B^{\top} B )^{- 1}$.
    For $i = 1, \dots , N$ and $p = 2, 3, \dots$, it holds that
    \begin{equation}  \label{diff_lml_p}
          \frac{d}{d\lambda} \sum_{i = 1}^{N} \frac{1}{( \lambda _{i} - \lambda )^{p - 1}}
        = ( p - 1 )          \sum_{i = 1}^{N} \frac{1}{( \lambda _{i} - \lambda )^{p    }}.
    \end{equation}
    This relationship implies that    we have   the summation $\sum_{i = 1}^{N} ( \lambda _{i} - \lambda )^{- p}$ for $p = 2, 3, \dots$ by differentiating (\ref{Eq_diff_1}) repeatedly.
    From this summation,              we obtain the summation $\sum_{i = 1}^{N}   \lambda _{i}            ^{- p}$ which is equal to the trace of $( B^{\top} B )^{- p}$
                         by substitution of $\lambda = 0$.
    To make handling of the right-hand-side of (\ref{Eq_diff_1}) easier,
    let us introduce functions $\hat{H}_{i}^{( p )}$ of $\lambda$ for $i = 1, \dots , N$ and $p = 1, 2, \dots$ defined by
    \begin{equation}  \label{defH}
        \hat{H}_{i}^{( p )} = \left \{ \begin{array}{ll}
                                  \displaystyle - \frac{  \hat{q}_{i}^{(     1 )}}{  \hat{q}_{i}^{( 0 )}},  &  p = 1,             \\
                                  \displaystyle   \frac{d \hat{H}_{i}^{( p - 1 )}}{d \lambda            },  &  p = 2, 3, \dots .
                              \end{array} \right .
    \end{equation}
    We can readily verify that it holds that
    \begin{equation}  \label{lambda_H}
          ( p - 1 )! \sum_{i = 1}^{N} \frac{1}{( \lambda _{i} - \lambda )^{p}}
        =            \sum_{i = 1}^{N} \hat{H}_{i}^{( p )}
    \end{equation}
    for $p = 1, 2, \dots$ by differentiating (\ref{Eq_diff_1}) repeatedly and taking care of (\ref{Eq_diff_1}), (\ref{diff_lml_p}) and (\ref{defH}).
    Thus, the trace of $( B^{\top} B )^{- p}$ is obtained by substituting $\lambda = 0$ into (\ref{lambda_H}).
    Let us introduce constants $H_{i}^{( p )}$ for $i = 1, \dots , N$ and $p = 1, 2, \dots$ defined by $H_{i}^{( p )} = \left. \hat{H}_{i}^{( p )} \right| _{\lambda = 0}$.
    The trace of $( B^{\top} B )^{- p}$ is expressed as
    \begin{equation}  \label{TrH_all}
        \textrm{Tr} ( ( B^{\top} B )^{- p} ) = \frac{1}{( p - 1 )!} \sum_{i = 1}^{N} H_{i}^{( p )},  \qquad  p = 1, 2, \dots .
    \end{equation}
    Thus, we can obtain these traces if the constants $H_{i}^{( p )}$ for $i = 1, \dots , N$ and $p = 1, 2, \dots$ are obtained by some means.
    The relationship (\ref{TrH_all}) implies that a formula for $H_{i}^{( p )}$ is required.
    We derive a recurrence relation for $\hat{H}_{i}^{( p )}$.
    A recurrence relation for $H_{i}^{( p )}$ is obtained by substitution of $\lambda = 0$ into the recurrence relation for $\hat{H}_{i}^{( p )}$.
    On the functions $\hat{H}_{i}^{( p )}$, the following lemma holds.
    \begin{lem}  \label{LemmaH}
        Let functions $\hat{h}_{i}^{( p )}$ of $\lambda$ for $i = 1, \dots , N$ and $p = 1, 2, \dots$ be defined by
        \begin{equation}  \label{def_hat_h}
            \hat{h}_{i}^{( p )} = - \frac{\hat{q}_{i}^{( p )}}{\hat{q}_{i}^{( 0 )}}.
        \end{equation}
        For $i = 1, \dots , N$, it holds
        \begin{align}
            &\hat{H}_{i}^{( 1 )} = \hat{h}_{i}^{( 1 )},                                                                                      & &                   \label{EqHh1}    \\
            &\hat{H}_{i}^{( p )} = \hat{h}_{i}^{( p )} + \sum_{k = 1}^{p - 1} {}_{p - 1} C_{k} \hat{h}_{i}^{( k )} \hat{H}_{i}^{( p - k )},  & &p = 2, 3, \dots .  \label{RR_hatG}
        \end{align}
    \end{lem}

    \noindent \textit{Proof}.

    In this proof, let $i = 1, \dots , N$.

    Substituting $p = 1$ into (\ref{def_hat_h}) and comparing with (\ref{defH}), we have (\ref{EqHh1}).

    We give the derivative of $\hat{h}_{i}^{( r )}$ for $r = 1, 2, \dots$.
    It holds
    \begin{equation*}
            \frac{d \hat{h}_{i}^{( r     )}                    }{d \lambda                  }
        = - \frac{  \hat{q}_{i}^{( r + 1 )}                    }{  \hat{q}_{i}^{( 0 )}      }
          + \frac{  \hat{q}_{i}^{( r     )} \hat{q}_{i}^{( 1 )}}{( \hat{q}_{i}^{( 0 )} )^{2}},  \qquad  r = 1, 2, \dots
    \end{equation*}
    from (\ref{def_hat_q_p}) and (\ref{def_hat_h}).
    Then, it holds that
    \begin{equation}  \label{deriv_hh_h}
        \frac{d \hat{h}_{i}^{( r )}}{d \lambda} = \hat{h}_{i}^{( r + 1 )} + \hat{h}_{i}^{( 1 )} \hat{h}_{i}^{( r )},  \qquad  r = 1, 2, \dots
    \end{equation}
    from (\ref{def_hat_h}).
    Using $\hat{H}_{i}^{( 1 )} = \hat{h}_{i}^{( 1 )}$ in (\ref{EqHh1}), we have another form
    \begin{equation}  \label{deriv_hh_H}
        \frac{d \hat{h}_{i}^{( r )}}{d \lambda} = \hat{h}_{i}^{( r + 1 )} + \hat{H}_{i}^{( 1 )} \hat{h}_{i}^{( r )},  \qquad  r = 1, 2, \dots .
    \end{equation}
    
    We use mathematical induction for proof.

    We write the definition of $\hat{H}_{i}^{( p )}$ for $p = 2, 3, \dots$ again.
    The definition is
    \begin{equation}  \label{defH_high_again}
        \hat{H}_{i}^{( p )} = \frac{d \hat{H}_{i}^{( p - 1 )}}{d \lambda}.
    \end{equation}
    
    We derive (\ref{RR_hatG}) for $p = 2$.
    Differentiating (\ref{EqHh1}) and using (\ref{deriv_hh_H}) and (\ref{defH_high_again}), we obtain
    \begin{equation}  \label{hatH_2_result}
        \hat{H}_{i}^{( 2 )} = \hat{h}_{i}^{( 2 )} + \hat{h}_{i}^{( 1 )} \hat{H}_{i}^{( 1 )}.
    \end{equation}
    Thus, (\ref{RR_hatG}) holds for $p = 2$.

    Hereafter, let $r$ be an integer such that $r \geq 2$ in this proof.
    Assume that (\ref{RR_hatG}) holds for $p = 2, \dots, r$.
    We consider differentiation of the function $\hat{H}_{i}^{( r )}$.
    Differentiating (\ref{RR_hatG}) for $p = r$ and using (\ref{deriv_hh_h}) and (\ref{defH_high_again}), we derive
    \begin{equation}  \label{diff_H_r}
        \hat{H}_{i}^{( r + 1 )} =                                                         \hat{h}_{i}^{( r + 1 )}
                                  +                                                       \hat{h}_{i}^{(     1 )}
                                                                                          \hat{h}_{i}^{( r     )}
                                  + \sum_{k = 1}^{r - 1} {}_{r - 1} C_{k} \left( \left(   \hat{h}_{i}^{( k + 1 )}
                                                                                        + \hat{h}_{i}^{(     1 )}
                                                                                          \hat{h}_{i}^{( k     )} \right) \hat{H}_{i}^{( r     - k )}
                                                                                        + \hat{h}_{i}^{( k     )}         \hat{H}_{i}^{( r + 1 - k )} \right) .
    \end{equation}
    Since it holds
    \begin{equation*}
          \hat{h}_{i}^{( 1 )} \hat{h}_{i}^{( r )} + \sum_{k = 1}^{r - 1} {}_{r - 1} C_{k} \hat{h}_{i}^{( 1 )} \hat{h}_{i}^{( k )} \hat{H}_{i}^{( r - k )}
        = \hat{h}_{i}^{( 1 )} \hat{H}_{i}^{( r )}
    \end{equation*}
    from the assumption, (\ref{diff_H_r}) is rewritten as
    \begin{equation}  \label{diff_H_r_ver1}
        \hat{H}_{i}^{( r + 1 )} =                                          \hat{h}_{i}^{( r + 1 )}
                                  +                                        \hat{h}_{i}^{(     1 )} \hat{H}_{i}^{( r         )}
                                  + \sum_{k = 1}^{r - 1} {}_{r - 1} C_{k}  \hat{h}_{i}^{( k + 1 )} \hat{H}_{i}^{( r     - k )}
                                  + \sum_{k = 1}^{r - 1} {}_{r - 1} C_{k}  \hat{h}_{i}^{( k     )} \hat{H}_{i}^{( r + 1 - k )}.
    \end{equation}
    We rearrange the third term in the right-hand-side of (\ref{diff_H_r_ver1}).
    Since it holds that
    \begin{equation*}
          \sum_{k          = 1}^{r - 1} {}_{r - 1} C_{k             } \hat{h}_{i}^{( k          + 1 )} \hat{H}_{i}^{( r     - k          )}
        = \sum_{k^{\prime} = 2}^{r    } {}_{r - 1} C_{k^{\prime} - 1} \hat{h}_{i}^{( k^{\prime}     )} \hat{H}_{i}^{( r + 1 - k^{\prime} )},
    \end{equation*}
    we have
    \begin{equation}  \label{slide_index}
                                                    \sum_{k = 1}^{r - 1} {}_{r - 1} C_{k    } \hat{h}_{i}^{( k + 1 )} \hat{H}_{i}^{( r     - k )}
        = \hat{h}_{i}^{( r )} \hat{H}_{i}^{( 1 )} + \sum_{k = 2}^{r - 1} {}_{r - 1} C_{k - 1} \hat{h}_{i}^{( k     )} \hat{H}_{i}^{( r + 1 - k )}.
    \end{equation}
    The fourth term in the right-hand-side of (\ref{diff_H_r_ver1}) is rewritten as
    \begin{equation}  \label{separate_terms}
                                                              \sum_{k = 1}^{r - 1} {}_{r - 1} C_{k}  \hat{h}_{i}^{( k     )} \hat{H}_{i}^{( r + 1 - k )}
        = ( r - 1 ) \hat{h}_{i}^{( 1 )} \hat{H}_{i}^{( r )} + \sum_{k = 2}^{r - 1} {}_{r - 1} C_{k}  \hat{h}_{i}^{( k     )} \hat{H}_{i}^{( r + 1 - k )}.
    \end{equation}
    From (\ref{diff_H_r_ver1}), (\ref{slide_index}) and (\ref{separate_terms}), we derive
    \begin{equation}  \label{diff_H_r_ver2}
        \hat{H}_{i}^{( r + 1 )} =                                                     \hat{h}_{i}^{( r + 1 )}
                                  +                                                   \hat{h}_{i}^{( r     )} \hat{H}_{i}^{(     1     )}
                                  +   \sum_{k = 2}^{r - 1} (   {}_{r - 1} C_{k - 1}
                                                             + {}_{r - 1} C_{k    } ) \hat{h}_{i}^{( k     )} \hat{H}_{i}^{( r + 1 - k )}
                                  + r                                                 \hat{h}_{i}^{(     1 )} \hat{H}_{i}^{( r         )}.
    \end{equation}
    It can readily be verified that the summation of the combinations in (\ref{diff_H_r_ver2}) is
    \begin{equation}  \label{RecurrCombin}
        {}_{r - 1} C_{k - 1} + {}_{r - 1} C_{k} = {}_{r} C_{k}
    \end{equation}
    in the case of $r \geq 3$.
    In the case of $r =    2$, the summation of the third term in the right-hand-side of (\ref{diff_H_r_ver2}) is zero.
    Then, we finally obtain
    \begin{equation*}
        \hat{H}_{i}^{( r + 1 )} = \hat{h}_{i}^{( r + 1 )} + \sum_{k = 1}^{r} {}_{r} C_{k} \hat{h}_{i}^{( k )} \hat{H}_{i}^{( r + 1 - k )}.
    \end{equation*}
    Thus, (\ref{RR_hatG}) holds for $p = r + 1$. ~ $\square$  \\

    By Lemma \ref{LemmaH}, we obtain a recurrence relation for $\hat{H}_{i}^{( p )}$.
    However, we have not obtained a recurrence relation for the functions $\hat{h}_{i}^{( p )}$.
    We show the following lemma which gives a method to compute $\hat{h}_{i}^{( p )}$ in a form of a recurrence relation.
    \begin{lem}  \label{Lemmah}
        The functions $\hat{h}_{i}^{( p )}$ for $i = 1, \dots , N$ and $p = 1, 2, \dots$ satisfy the following recurrence relation.
        For $p = 1$, the recurrence relation is
        \begin{align}
            &\hat{h}_{1}^{( 1 )} =   \displaystyle \frac{1                                                  }{\hat{q}_{1}^{( 0 )}},  & &                   \label{h1_1}  \\
            &\hat{h}_{i}^{( 1 )} =   \displaystyle \frac{\hat{e}_{i - 1}^{( 0 )} \hat{h}_{i - 1}^{( 1 )} + 1}{\hat{q}_{i}^{( 0 )}},  & &i = 2, \dots , N.  \label{h1_i}
        \end{align}
        For $i = 1$ and $p = 2, 3, \dots$, the recurrence relation is
        \begin{equation}  \label{h1p_0}
            \hat{h}_{1}^{( p )} = 0.
        \end{equation}
        For $i = 2, \dots , N$ and $p = 2, 3, \dots$, the recurrence relation is
        \begin{equation}  \label{hi_high}
            \hat{h}_{i}^{( p )} =   \frac{\hat{e}_{i - 1}^{( 0 )}}{\hat{q}_{i}^{( 0 )}} \left(     \hat{h}_{i - 1}^{( p     )}
                                                                                               + p \hat{h}_{i - 1}^{(     1 )}
                                                                                                   \hat{h}_{i - 1}^{( p - 1 )} \right)
                                  + \sum_{k = 1}^{p - 2} {}_{p} C_{k}                              \hat{h}_{i - 1}^{( k     )}
                                                                                                   \hat{h}_{i    }^{( p - k )}.
        \end{equation}
    \end{lem}

    \noindent \textit{Proof}.

    In this proof, another key point of the derivation of the recurrence relation is applied.
    The key point is as follows.
    We differentiate the relationships     $q_{i} + e_{i - 1} - \lambda = \hat{q}_{i}^{( 0 )} + \hat{e}_{i - 1}^{( 0 )}$ for $i = 1, \dots , N$     shown in (\ref{eq_hat_sum})
                                       and $q_{i}   e_{i    }           = \hat{q}_{i}^{( 0 )}   \hat{e}_{i    }^{( 0 )}$ for $i = 1, \dots , N - 1$ shown in (\ref{eq_hat_prod}).
    Then, we derive a recurrence relation which the functions $\hat{h}_{i}^{( p )}$ satisfy.

    We introduce functions $\hat{e}_{i}^{( p )}$ of $\lambda$ for $i = 0, 1, \dots , N - 1$ and $p = 1, 2, \dots$ defined by
    \begin{equation*}
        \hat{e}_{i}^{( p )} = \frac{d^{p} \hat{e}_{i}^{( 0 )}}{d \lambda ^{p}}.
    \end{equation*}
    Differentiating (\ref{eq_hat_prod}) repeatedly, we have
    \begin{equation*}
        \sum_{k = 0}^{p} {}_{p} C_{k} \hat{q}_{i}^{( k )} \hat{e}_{i}^{( p - k )} = 0,  \qquad  i = 1, \dots N - 1,  \qquad  p = 1, 2, \dots .
    \end{equation*}
    Solving this equation for $\hat{e}_{i}^{( p )}$, we obtain
    \begin{equation}  \label{diff_prod_hat}
        \hat{e}_{i}^{( p )} = \sum_{k = 1}^{p} {}_{p} C_{k} \hat{h}_{i}^{( k )} \hat{e}_{i}^{( p - k )},  \qquad  i = 1, \dots N - 1,  \qquad  p = 1, 2, \dots
    \end{equation}
    since $\hat{h}_{i}^{( k )} = - \hat{q}_{i}^{( k )} / \hat{q}_{i}^{( 0 )}$ from the definition.
    
    We show that (\ref{h1_1}) and (\ref{h1_i}) hold.
    Differentiating (\ref{eq_hat_sum}), we have
    \begin{equation}  \label{diff_qe_1}
        \hat{q}_{i}^{( 1 )} + \hat{e}_{i - 1}^{( 1 )} = - 1,  \qquad  i = 1, \dots , N.
    \end{equation}
    Then, $\hat{h}_{i}^{( 1 )} = - \hat{q}_{i}^{( 1 )} / \hat{q}_{i}^{( 0 )}$ is expressed as
    \begin{equation}  \label{h1_qe}
        \hat{h}_{i}^{( 1 )} = \frac{\hat{e}_{i - 1}^{( 1 )} + 1}{\hat{q}_{i}^{( 0 )}},  \qquad  i = 1, \dots , N.
    \end{equation}
    Substituting $p = 1$ into (\ref{diff_prod_hat}), we have
    \begin{equation}  \label{e1}
        \hat{e}_{i}^{( 1 )} = \hat{h}_{i}^{( 1 )} \hat{e}_{i}^{( 0 )},  \qquad  i = 1, \dots , N - 1.
    \end{equation}
    From (\ref{h1_qe}), (\ref{e1}) and $\hat{e}_{0}^{(0)} = 0$, we immediately obtain (\ref{h1_1}) and (\ref{h1_i}).

    Hereafter, let $p = 2, 3 , \dots$ in this proof.
    We show that the relationships (\ref{h1p_0}) and (\ref{hi_high}) hold.
    Differentiating (\ref{diff_qe_1}) repeatedly, we obtain
    \begin{equation}  \label{diff_qpe_p}
        \hat{q}_{i}^{( p )} + \hat{e}_{i - 1}^{( p )} = 0,  \qquad  i = 1, \dots , N.
    \end{equation}
    From this relationship and the definition $\hat{h}_{i}^{( p )} = - \hat{q}_{i}^{( p )} / \hat{q}_{i}^{( 0 )}$, the functions $\hat{h}_{i}^{( p )}$ for $p = 2, 3, \dots$ are
    \begin{equation}  \label{h_ratio_qe}
        \hat{h}_{i}^{( p )} = \frac{\hat{e}_{i - 1}^{( p )}}{\hat{q}_{i}^{( 0 )}},  \qquad  i = 1, \dots , N.
    \end{equation}

    We show that (\ref{h1p_0}) holds.
    It holds that $\hat{e}_{0}^{( p )} = 0$ since $\hat{e}_{0}^{( 0 )} = 0$.
    Then, substituting $i = 1$ into (\ref{h_ratio_qe}), we have (\ref{h1p_0}).

    We show that (\ref{hi_high}) holds.
    Hereafter, let $i = 2, \dots N$ in this proof.
    Substituting (\ref{diff_prod_hat}) into (\ref{h_ratio_qe}), we obtain
    \begin{equation}  \label{hp_sum_ek}
        \hat{h}_{i}^{( p )} = \frac{1}{\hat{q}_{i}^{( 0 )}} \sum_{k = 1}^{p} {}_{p} C_{k} \hat{h}_{i - 1}^{( k )} \hat{e}_{i - 1}^{( p - k )}.
    \end{equation}
    It follows from (\ref{e1}) that
    \begin{equation}  \label{rearr_1_pm1}
                                    {}_{p} C_{p - 1} \hat{h}_{i - 1}^{( p - 1 )} \hat{e}_{i - 1}^{( 1 )}
        = p \hat{e}_{i - 1}^{( 0 )}                  \hat{h}_{i - 1}^{( p - 1 )} \hat{h}_{i - 1}^{( 1 )}.
    \end{equation}
    From (\ref{hp_sum_ek}) and (\ref{rearr_1_pm1}), it holds
    \begin{equation}  \label{hip_1st}
                   \hat{h}_{i}^{( p )}
        = \frac{1}{\hat{q}_{i}^{( 0 )}}
          \left(                                                                                           \hat{h}_{i - 1}^{( p     )} \hat{e}_{i - 1}^{( 0     )}
                 + p                                   \hat{e}_{i - 1}^{( 0 )} \hat{h}_{i - 1}^{( p - 1 )} \hat{h}_{i - 1}^{( 1     )} 
                 +   \sum_{k = 1}^{p - 2} {}_{p} C_{k}                         \hat{h}_{i - 1}^{( k     )}                             \hat{e}_{i - 1}^{( p - k )} \right) .
    \end{equation}
    Then, we obtain (\ref{hi_high}) from (\ref{h_ratio_qe}) and (\ref{hip_1st}). ~ $\square$  \\
    
    \begin{rmk} \label{Rmk_Hh1}
        From Lemmas \ref{LemmaH} and \ref{Lemmah}, the functions $\hat{H}_{i}^{( 1 )}$ satisfy the following recurrence relation
        \begin{align*}
            &\hat{H}_{1}^{( 1 )} = \displaystyle \frac{1                                                  }{\hat{q}_{1}^{( 0 )}},  & &                   \\
            &\hat{H}_{i}^{( 1 )} = \displaystyle \frac{\hat{e}_{i - 1}^{( 0 )} \hat{H}_{i - 1}^{( 1 )} + 1}{\hat{q}_{i}^{( 0 )}},  & &i = 2, \dots , N.
        \end{align*}
    \end{rmk}
    ~\\

    Obviously, it holds that     $q_{i} = \left. \hat{q}_{i}^{( 0 )} \right| _{\lambda = 0}$ for $i = 1, \dots , N$
                             and $e_{i} = \left. \hat{e}_{i}^{( 0 )} \right| _{\lambda = 0}$ for $i = 1, \dots , N - 1$.
    Substituting $\lambda = 0$ into the recurrence relations in Lemmas \ref{LemmaH} and \ref{Lemmah} and Remark \ref{Rmk_Hh1} and considering (\ref{TrH_all}),
    we finally obtain one of the main theorems in this paper.

    \begin{thm} \label{Th_G}
        Let $B$ be an upper bidiagonal matrix defined in (\ref{UpperBidiagonal}).

        Let us introduce constants $\tilde{F}_{i}$ for $i = 2, \dots , N$ defined as
        \begin{equation*}
            \tilde{F}_{i} = \frac{e_{i - 1}}{q_{i}},  \qquad  i = 2, \dots , N.
        \end{equation*}
        For $i = 1, \dots , N$ and $p = 1, 2, \dots$, let $h_{i}^{( p )}$ be constants which satisfy the following recurrence relation.
        For $p = 1$, the recurrence relation is
        \begin{gather*}
            h_{1}^{( 1 )} =                                   \frac{1}{q_{1}},                             \\
            h_{i}^{( 1 )} = \tilde{F}_{i} h_{i - 1}^{( 1 )} + \frac{1}{q_{i}},  \qquad  i = 2, \dots , N.
        \end{gather*}
        For $i = 1$ and $p = 2, 3, \dots$, the recurrence relation is
        \begin{equation*}  \label{Th_h1p_0}
            h_{1}^{( p )} =  0.
        \end{equation*}
        For $i = 2, \dots , N$ and $p = 2, 3, \dots$, the recurrence relation is
        \begin{equation*}
            h_{i}^{( p )} = \tilde{F}_{i} \left(                                       h_{i - 1}^{( p     )}
                                                 + p                                   h_{i - 1}^{(     1 )}
                                                                                       h_{i - 1}^{( p - 1 )} \right)
                                                 +   \sum_{k = 1}^{p - 2} {}_{p} C_{k} h_{i - 1}^{( k     )}
                                                                                       h_{i    }^{( p - k )}.
        \end{equation*}

        For $i = 1, \dots , N$, let $H_{i}^{( 1 )}$ be constants given as
        \begin{equation*}
            H_{i}^{( 1 )} = h_{i}^{( 1 )}.
        \end{equation*}
        For $i = 1, \dots , N$ and $p = 2, 3, \dots$, let $H_{i}^{( p )}$ be constants which satisfy the following recurrence relation
        \begin{equation*}
            H_{i}^{( p )} = h_{i}^{( p )} + \sum_{k = 1}^{p - 1} {}_{p - 1} C_{k} h_{i}^{( k )} H_{i}^{( p - k )}.
        \end{equation*}

        The traces $\textrm{Tr} ( ( B^{\top} B )^{- p} )$ for $p = 1, 2, \dots$ are computed by
        \begin{equation*}
            \textrm{Tr} ( ( B^{\top} B )^{- p} ) = \frac{1}{( p - 1 ) !} \sum_{i = 1}^{N} H_{i}^{( p )}.
        \end{equation*}
    \end{thm}

    The formula in Theorem \ref{Th_G} consists of only summation, multiplication and division among positive quantities.
    Then, possibility of cancellation error is clearly excluded.
    \begin{rmk}
        The recurrence relation of $h_{i}^{( 1 )}$ $( i = 1, \dots , N )$ in Theorem \ref{Th_G} is equivalent
        to the recurrence relation for the diagonal entries of $( B B^{\top} )^{- 1}$ shown in Remark 4.6 in \cite{KYN11}.
        Then, the constants $h_{i}^{( 1 )}$ and $H_{i}^{( 1 )}$ are the $( i, i )$-entry of $( B B^{\top} )^{- 1}$.
        See also Remark 4.7 in \cite{KYN11}.
    \end{rmk}

    \subsection{Derivation - type II}  \label{Sect_DerivII}
    
    In this subsection, we consider the matrix $B B^{\top}$ instead of the matrix $B^{\top} B$.
    A recurrence relation for computation of the traces $\textrm{Tr} ( ( B B^{\top} )^{- p} )$ for an arbitrary positive integer $p$ is derived.
    Note that it holds $\textrm{Tr} ( ( B^{\top} B )^{- p} ) = \textrm{Tr} ( ( B B^{\top} )^{- p} )$.
    Let the eigenvalues of $B B^{\top}$ be $\tilde{\lambda}_{1}, \dots , \tilde{\lambda}_{N}$.
    For an arbitrary positive integer $p$, the eigenvalues of $( B B^{\top} )^{- p}$ are $\tilde{\lambda}_{1}^{- p}, \dots , \tilde{\lambda}_{N}^{- p}$.
    Then, the summation $\sum_{i = 1}^{N} \tilde{\lambda}_{i}^{- p}$ is  the trace of $( B B^{\top} )^{- p}$.
    We derive a formula to compute this summation.
    The procedure of the derivation in this subsection is similar to that in the previous subsection.
    Then, we show only an outline of the derivation.
    The matrix $B B^{\top}$ is given as
    \begin{equation*}
        B B^{\top} = \left( \begin{array}{cccc}
                                      q_{1} + e_{1}   &  \sqrt{q_{2} e_{1}}  &                                &                          \\
                                \sqrt{q_{2}   e_{1}}  &  \ddots              &  \ddots                        &                          \\
                                                      &  \ddots              &        q_{N - 1} + e_{N - 1}   &  \sqrt{q_{N} e_{N - 1}}  \\
                                                      &                      &  \sqrt{q_{N    }   e_{N - 1}}  &        q_{N}
                            \end{array}                                                                                                      \right) .
    \end{equation*}
    It can be readily verified that we obtain the following matrix
    \begin{equation*}
        \tilde{A} = \left( \begin{array}{cccc}
                               q_{1} + e_{1}  &  q_{2} e_{1}  &                         &                   \\
                               1              &  \ddots       &  \ddots                 &                   \\
                                              &  \ddots       &  q_{N - 1} + e_{N - 1}  &  q_{N} e_{N - 1}  \\
                                              &               &  1                      &  q_{N}
                           \end{array}                                                                          \right)
    \end{equation*}
    by similarity transformation.
    Let us consider decomposition of the matrix
    \begin{equation*}
        \tilde{A} - \lambda I = \left( \begin{array}{cccc}
                                           q_{1} + e_{1} - \lambda  &  q_{2} e_{1}  &                                   &                             \\
                                           1                        &  \ddots       &  \ddots                           &                             \\
                                                                    &  \ddots       &  q_{N - 1} + e_{N - 1} - \lambda  &  q_{N} e_{N - 1}            \\
                                                                    &               &  1                                &  q_{N}           - \lambda
                                       \end{array}                                                                                                        \right)
    \end{equation*}
    into the matrix product expressed as
    \begin{equation*}
        \tilde{A} - \lambda I = \left( \begin{array}{cccc}
                                           1  &  \check{e}_{1}^{( 0 )}  &          &                             \\
                                              &  1                      &  \ddots  &                             \\
                                              &                         &  \ddots  &  \check{e}_{N - 1}^{( 0 )}  \\
                                              &                         &          &  1
                                       \end{array}                                                                            \right)
                                \left( \begin{array}{cccc}
                                            \check{q}_{1}^{( 0 )}  &                         &          &                         \\
                                            1                      &  \check{q}_{2}^{( 0 )}  &          &                         \\
                                                                   &  \ddots                 &  \ddots  &                         \\
                                                                   &                         &  1       &  \check{q}_{N}^{( 0 )}
                                       \end{array}                                                                                   \right) ,
    \end{equation*}
    where $\check{q}_{i}^{( 0 )}$ for $i = 1, \dots , N$ and $\check{e}_{i}^{( 0 )}$ for $i = 1, \dots , N - 1$ are functions of $\lambda$.

    Similarly to the previous subsection, we introduce functions $\check{q}_{i}^{( p )}$ and $\check{H}_{i}^{( p )}$ of $\lambda$ for $i = 1, \dots , N$ and $p = 1, 2, \dots$.
    The definition of $\check{q}_{i}^{( p )}$ is
    \begin{equation*}
        \check{q}_{i}^{( p )} = \frac{d^{p} \check{q}_{i}^{( 0 )}}{d \lambda ^{p}}.
    \end{equation*}
    The definition of $\check{H}_{i}^{( p )}$ is
    \begin{equation*}
        \check{H}_{i}^{( p )} = \left \{ \begin{array}{ll}
                                             \displaystyle - \frac{  \check{q}_{i}^{(     1 )}}{  \check{q}_{i}^{( 0 )}},  &  p = 1,             \\
                                             \displaystyle   \frac{d \check{H}_{i}^{( p - 1 )}}{d \lambda              },  &  p = 2, 3, \dots .
                                         \end{array} \right .
    \end{equation*}
    Similarly to the derivation of (\ref{lambda_H}) in the previous subsection, we have
    \begin{equation*}
          ( p - 1 )! \sum_{i = 1}^{N} \frac{1}{( \tilde{\lambda}_{i} - \lambda )^{p}}
        =            \sum_{i = 1}^{N} \check{H}_{i}^{( p )},                           \qquad  p = 1, 2, \dots .
    \end{equation*}
    Let us introduce constants $\tilde{H}_{i}^{( p )}$ for $i = 1, \dots , N$ and $p = 1, 2, \dots$
    defined by $\tilde{H}_{i}^{( p )} = \left. \check{H}_{i}^{( p )} \right| _{\lambda = 0}$.
    The trace of $( B B^{\top} )^{- p}$ is expressed as
    \begin{equation}  \label{TrtH_all}
        \textrm{Tr} ( ( B B^{\top} )^{- p} ) = \frac{1}{( p - 1 )!} \sum_{i = 1}^{N} \tilde{H}_{i}^{( p )},  \qquad  p = 1, 2, \dots .
    \end{equation}
    We derive a recurrence relation for $\check{H}_{i}^{( p )}$.
    A recurrence relation for $\tilde{H}_{i}^{( p )}$ is obtained by substitution of $\lambda = 0$ into the recurrence relation for $\check{H}_{i}^{( p )}$.
    The following lemma holds.
    \begin{lem}  \label{LemmatH}
        Let functions $\check{h}_{i}^{( p )}$ of $\lambda$ for $i = 1, \dots , N$ and $p = 1, 2, \dots$ be defined by
        \begin{equation*}
            \check{h}_{i}^{( p )} = - \frac{\check{q}_{i}^{( p )}}{\check{q}_{i}^{( 0 )}}.
        \end{equation*}
        For $i = 1, \dots , N$, it holds
        \begin{align*}
            &\check{H}_{i}^{( 1 )} = \check{h}_{i}^{( 1 )},                                                                                                                    \\
            &\check{H}_{i}^{( p )} = \check{h}_{i}^{( p )} + \sum_{k = 1}^{p - 1} {}_{p - 1} C_{k} \check{h}_{i}^{( k )} \check{H}_{i}^{( p - k )},  \qquad  p = 2, 3, \dots .
        \end{align*}
    \end{lem}
    ~\\ 
    \noindent
    Proof of this lemma is similar to that of Lemma \ref{LemmaH}.

    We show the following lemma which gives a method to compute $\check{h}_{i}^{( p )}$ in a form of recurrence relation.
    \begin{lem}  \label{Lemmath}
        The functions $\check{h}_{i}^{( p )}$ for $i = 1, \dots , N$ and $p = 1, 2, \dots$ satisfy the following recurrence relation.
        For $p = 1$, the recurrence relation is
        \begin{align*}
            &\check{h}_{N}^{( 1 )} =   \displaystyle \frac{1                                                  }{\check{q}_{N}^{( 0 )}},  & &                       \\
            &\check{h}_{i}^{( 1 )} =   \displaystyle \frac{\check{e}_{i}^{( 0 )} \check{h}_{i + 1}^{( 1 )} + 1}{\check{q}_{i}^{( 0 )}},  & &i = 1, \dots , N - 1.
        \end{align*}
        For $i = N$ and $p = 2, 3, \dots$, the recurrence relation is
        \begin{equation*}
            \check{h}_{N}^{( p )} = 0.
        \end{equation*}
        For $i = 1, \dots , N - 1$ and $p = 2, 3, \dots$, the recurrence relation is
        \begin{equation*}
            \check{h}_{i}^{( p )} =   \frac{\check{e}_{i}^{( 0 )}}{\check{q}_{i}^{( 0 )}} \left(     \check{h}_{i + 1}^{( p     )}
                                                                                                 + p \check{h}_{i + 1}^{(     1 )}
                                                                                                     \check{h}_{i + 1}^{( p - 1 )} \right)
                                    + \sum_{k = 1}^{p - 2} {}_{p} C_{k}                              \check{h}_{i + 1}^{( k     )}
                                                                                                     \check{h}_{i    }^{( p - k )}.
        \end{equation*}
    \end{lem}
    Proof of this lemma is similar to that of Lemma \ref{Lemmah}.

    \begin{rmk} \label{Rmk_tHh1}
        From Lemmas \ref{LemmatH} and \ref{Lemmath}, the functions $\check{H}_{i}^{( 1 )}$ satisfy the following recurrence relation
        \begin{align*}
            &\check{H}_{N}^{( 1 )} = \displaystyle \frac{1                                                  }{\check{q}_{N}^{( 0 )}},  & &                       \\
            &\check{H}_{i}^{( 1 )} = \displaystyle \frac{\check{e}_{i}^{( 0 )} \check{H}_{i + 1}^{( 1 )} + 1}{\check{q}_{i}^{( 0 )}},  & &i = 1, \dots , N - 1.
        \end{align*}
    \end{rmk}
    ~\\ 

    Substituting $\lambda = 0$ into the recurrence relation in Lemmas \ref{LemmatH} and \ref{Lemmath} and Remark \ref{Rmk_tHh1} and considering (\ref{TrtH_all}), 
    we finally obtain one of the main theorems in this paper.

    \begin{thm} \label{Th_tG}
        Let $B$ be an upper bidiagonal matrix defined in (\ref{UpperBidiagonal}).

        Let us introduce constants $F_{i}$ for $i = 1, \dots , N - 1$ defined as
        \begin{equation*}
            F_{i} = \frac{e_{i}}{q_{i}},  \qquad  i = 1, \dots , N - 1.
        \end{equation*}
        For $i = 1, \dots , N$ and $p = 1, 2, \dots$, let $\tilde{h}_{i}^{( p )}$ be constants which satisfy the following recurrence relation.
        For $p = 1$, the recurrence relation is
        \begin{gather*}
            \tilde{h}_{N}^{( 1 )} =                                   \frac{1}{q_{N}},                                 \\
            \tilde{h}_{i}^{( 1 )} = F_{i} \tilde{h}_{i + 1}^{( 1 )} + \frac{1}{q_{i}},  \qquad  i = 1, \dots , N - 1.
        \end{gather*}
        For $i = N$ and $p = 2, 3, \dots$, the recurrence relation is
        \begin{equation*}
            \tilde{h}_{N}^{( p )} =  0.
        \end{equation*}
        For $i = 1, \dots , N - 1$ and $p = 2, 3, \dots$, the recurrence relation is
        \begin{equation*}
            \tilde{h}_{i}^{( p )} = F_{i} \left(                                       \tilde{h}_{i + 1}^{( p     )}
                                                 + p                                   \tilde{h}_{i + 1}^{(     1 )}
                                                                                       \tilde{h}_{i + 1}^{( p - 1 )} \right)
                                                 +   \sum_{k = 1}^{p - 2} {}_{p} C_{k} \tilde{h}_{i + 1}^{( k     )}
                                                                                       \tilde{h}_{i    }^{( p - k )}.
        \end{equation*}

        For $i = 1, \dots , N$, let $\tilde{H}_{i}^{( 1 )}$ be constants given as
        \begin{equation*}
            \tilde{H}_{i}^{( 1 )} = \tilde{h}_{i}^{( 1 )}.
        \end{equation*}
        For $i = 1, \dots , N$ and $p = 2, 3, \dots$, let $\tilde{H}_{i}^{( p )}$ be constants which satisfy the following recurrence relation
        \begin{equation*}
            \tilde{H}_{i}^{( p )} = \tilde{h}_{i}^{( p )} + \sum_{k = 1}^{p - 1} {}_{p - 1} C_{k} \tilde{h}_{i}^{( k )} \tilde{H}_{i}^{( p - k )}.
        \end{equation*}

        The traces $\textrm{Tr} ( ( B B^{\top} )^{- p} )$ for $p = 1, 2, \dots$ are computed by
        \begin{equation*}
            \textrm{Tr} ( ( B B^{\top} )^{- p} ) = \frac{1}{( p - 1 ) !} \sum_{i = 1}^{N} \tilde{H}_{i}^{( p )}.
        \end{equation*}
    \end{thm}

    The formula in Theorem \ref{Th_tG} has the same merit as that the formula in Theorem \ref{Th_G} has.
    \begin{rmk}
        The recurrence relation of $\tilde{h}_{i}^{( 1 )}$ $( i = 1, \dots , N )$ in Theorem \ref{Th_tG} is equivalent
        to the recurrence relation for the diagonal entries of $( B^{\top} B )^{- 1}$ shown in Remark 4.6 in \cite{KYN11}.
        Then, the constants $\tilde{h}_{i}^{( 1 )}$ and $\tilde{H}_{i}^{( 1 )}$ are the $( i, i )$-entry of $( B^{\top} B )^{- 1}$.
    \end{rmk}

    \section{Computational costs for the traces}  \label{Cost}

    In this section, we discuss computational cost for the trace $\textrm{Tr} ( ( B^{\top} B )^{- M} )$.
    We consider the case where the matrix size $N$ of $B$ and the order $M$ are sufficiently large.
    In this section, let input be the diagonal and the upper subdiagonal entries of $B$.
    For $i = 1, \dots , N   $ , let $Q_{i}$ be $Q_{i} = \sqrt{q_{i}}$.
    For $i = 1, \dots , N- 1$ , let $E_{i}$ be $E_{i} = \sqrt{e_{i}}$.
    Let us introduce constants $\check{B}_{i}$ for $i = 1, \dots , N$ defined as
    \begin{equation}  \label{def_cB}
        \check{B}_{i} = q_{i}^{- 1}.
    \end{equation}
    An algorithm for computation of the trace $\textrm{Tr} ( ( B^{\top} B )^{- M} )$ is given in Algorithm \ref{Algo_Higher}.
    The lines from 1 to 9 compute the constants $\check{B}_{i}$, $h_{i}^{( 1 )}$ and $H_{i}^{( 1 )}$ for $i = 1, \dots , N$ and  $\tilde{F}_{i}$ for $i = 1, \dots , N - 1$.
    The constant $S_{k}^{( p )}$ represents ${}_{p} C_{k}$ for each $p$ and $k$.
    In the lines from 12 to 16, the constants ${}_{p} C_{k}$ for $k = 1, \dots , p$ are set.
    In the line 14, the relationship shown in (\ref{RecurrCombin}) is used.
    The lines from 17 to 30 compute the constants $h_{i}^{( p )}$ and $H_{i}^{( p )}$ for $i = 1, \dots , N$ and $p = 2, \dots , M$.
    In the lines from 20 to 23, the summation $\sum_{k = 1}^{p - 2} {}_{p    } C_{k} h_{i - 1}^{( k )} h_{i}^{( p - k )}$ is computed and is stored in the variable $tmp$.
    In the line             24, the constant $h_{i}^{( p )}$ is obtained.
    In the lines from 25 to 28, the summation $\sum_{k = 1}^{p - 1} {}_{p - 1} C_{k} h_{i    }^{( k )} H_{i}^{( p - k )}$ is computed and is stored in the variable $tmp$.
    In the line             29, the constant $H_{i}^{( p )}$ is obtained.
    The    lines from 32 to 35 are used to compute the trace.
    The variable $J$ is used to compute the trace.
    \begin{algorithm}
        \caption{Computation of the trace $\textrm{Tr} ( ( B^{\top} B )^{- M} )$ for sufficiently large $M$ and matrix size} \label{Algo_Higher}
        \begin{algorithmic}[1]
            \STATE $\check{B}_{1}         \gets 1.0                   / ( Q_{1} * Q_{1} )$
            \STATE $       h _{1}^{( 1 )} \gets \check{B}_{1}                            $
            \STATE $       H _{1}^{( 1 )} \gets        h _{1}^{( 1 )}                    $
            \FOR{$i = 2$ \TO $N$ {\bf by} $+ 1$}
                \STATE $\check{B}_{i}         \gets 1.0               / ( Q_{i    }         *        Q _{i} )$
                \STATE $\tilde{F}_{i}         \gets        E _{i - 1} *   E_{i - 1}         * \check{B}_{i}  $
                \STATE $       h _{i}^{( 1 )} \gets \tilde{F}_{i    } *   h_{i - 1}^{( 1 )} + \check{B}_{i}  $
                \STATE $       H _{i}^{( 1 )} \gets                       h_{i    }^{( 1 )}                  $
            \ENDFOR
            \FOR{$p = 2$ \TO $M$ {\bf by} $+ 1$}
                \STATE $p_{d} \gets p    $ : $p$ is cast into double precision number
                \STATE $S_{1}^{( p )} \gets p_{d}$
                \FOR{$k = 2$ \TO $p - 1$ {\bf by} $+ 1$}
                    \STATE $S_{k}^{( p )} \gets S_{k - 1}^{( p - 1 )} + S_{k}^{( p - 1 )}$
                \ENDFOR
                \STATE $S_{p}^{( p )} \gets 1.0$
                \STATE $h_{1}^{( p )} \gets 0.0$
                \STATE $H_{1}^{( p )} \gets S_{1}^{( p - 1 )} * h_{1}^{( 1 )} * H_{1}^{( p - 1 )}$
                \FOR{$i = 2$ \TO $N$ {\bf by} $+ 1$}
                    \STATE $tmp \gets 0.0$
                    \FOR{$k = 1$ \TO $p - 2$ {\bf by} $+ 1$}
                        \STATE $tmp \gets tmp + S_{k}^{( p )} * h_{i - 1}^{( k )} * h_{i}^{( p - k )}$
                    \ENDFOR
                    \STATE $h_{i}^{( p )} \gets \tilde{F}_{i} * ( h_{i - 1}^{( p )} + p_{d} * h_{i - 1}^{( 1 )} * h_{i - 1}^{( p - 1 )} ) + tmp$
                    \STATE $tmp \gets 0.0$
                    \FOR{$k = 1$ \TO $p - 1$ {\bf by} $+ 1$}
                        \STATE $tmp \gets tmp + S_{k}^{( p - 1 )} * h_{i}^{( k )} * H_{i}^{( p - k )}$
                    \ENDFOR
                    \STATE $H_{i}^{( p )} \gets h_{i}^{( p )} + tmp$
                \ENDFOR
            \ENDFOR
            \STATE $J \gets H_{1}^{( p )}$
            \FOR{$i = 2$ \TO $N$ {\bf by} $+ 1$}
                \STATE $J \gets J + H_{i}^{( p )}$
            \ENDFOR
            \RETURN $J$
        \end{algorithmic}
    \end{algorithm}

    Thus, the following remark follows.
    
    \begin{rmk}
        In Algorithm \ref{Algo_Higher}, there exist three nested loops.
        The first  loop  is  from the line 10 to the line 31.
        The second loop  is  from the line 19 to the line 30.
        The third  loops are from the line 21 to the line 23 and from the line 26 to the line 28.
        Then, the order of computational cost for the trace $\textrm{Tr} ( ( B^{\top} B )^{- M} )$ is $O ( M^{2} N )$.
    \end{rmk}

    \section{An efficient implementation of the formula for the trace $\textrm{Tr} ( ( B B^{\top} )^{- 2} )$}  \label{Implementation}

    In this section, we present an efficient implementation of the formula for the trace $\textrm{Tr} ( ( B B^{\top} )^{- 2} )$.
    Computation of this trace is practically important.
    For example, see \cite{vonMatt97, YKTN13}.

    It is obvious that a way for computing of such a trace with a smaller number of arithmetic operations is more desirable.
    Then, we rearrange the recurrence relation in Theorem \ref{Th_G} to reduce the number of arithmetic operations.
    Note that $H_{i}^{( 1 )} = h_{i}^{( 1 )}$ for $i = 1, \dots , N$ as defined in Theorem \ref{Th_G}.
    Then, the constants $H_{i}^{( 1 )}$ satisfy the following recurrence relation
    \begin{align*}
        &H_{1}^{( 1 )} =                                   \check{B}_{1},  & &                   \\
        &H_{i}^{( 1 )} = \tilde{F}_{i} H_{i - 1}^{( 1 )} + \check{B}_{i},  & &i = 2, \dots , N,
    \end{align*}
    where $\check{B}_{i}$ is defined in (\ref{def_cB}), from Theorem \ref{Th_G}.
    We rearrange the recurrence relation for $H_{i}^{( 2 )}$ and $h_{i}^{( 2 )}$ which are shown in Theorem \ref{Th_G}.
    The recurrence relation for $H_{i}^{( 2 )}$ is
    \begin{equation}  \label{H2_original}
        H_{i}^{( 2 )} = h_{i}^{( 2 )} + h_{i}^{( 1 )} H_{i}^{( 1 )},  \qquad  i = 1, \dots , N.
    \end{equation}
    The recurrence relation for $h_{i}^{( 2 )}$ is
    \begin{align}
        &h_{1}^{( 2 )} = 0,                                                                                          & &                   \label{h2_1}     \\
        &h_{i}^{( 2 )} = \tilde{F}_{i} \left( h_{i - 1}^{( 2 )} + 2 \left( h_{i - 1}^{( 1 )} \right) ^{2} \right) ,  & &i = 2, \dots , N.  \label{h2_2toN}
    \end{align}
    Using $H_{i}^{( 1 )} = h_{i}^{( 1 )}$, we readily derive
    \begin{equation}  \label{h2_simplified}
        h_{i}^{( 2 )} = \tilde{F}_{i} \left( H_{i - 1}^{( 2 )} + \left( H_{i - 1}^{( 1 )} \right) ^{2} \right) ,  \qquad  i = 2, \dots , N
    \end{equation}
    from (\ref{H2_original}) and (\ref{h2_2toN}).
    Substituting (\ref{h2_1}) or (\ref{h2_simplified}) into (\ref{H2_original}) and using $H_{i}^{( 1 )} = h_{i}^{( 1 )}$, we obtain
    \begin{align*}
        &H_{1}^{( 2 )} =                                          \left( H_{1    }^{( 1 )} \right) ^{2},                               \\
        &H_{i}^{( 2 )} = \tilde{F}_{i} \left( H_{i - 1}^{( 2 )} + \left( H_{i - 1}^{( 1 )} \right) ^{2} \right)
                                                                + \left( H_{i    }^{( 1 )} \right) ^{2},         & &i = 2, \dots , N.
    \end{align*}
    Then, the constants $h_{i}^{( 1 )}$ and $h_{i}^{( 2 )}$ are not necessary in computation of $H_{i}^{( 2 )}$.
    Let us introduce auxiliary constants $\Phi _{i}$ for $i = 1, \dots , N$ defined as
    \begin{equation}  \label{def_Phi}
        \Phi _{i} = \left( H_{i}^{( 1 )} \right) ^{2} ,  \qquad  i = 1, \dots , N.
    \end{equation}
    We have the following corollary of Theorem \ref{Th_G}.
    \begin{cor}  \label{Cor_G2}
        The constants $H_{i}^{( 2 )}$ for $i = 1, \dots , N$ are obtained from the following recurrence relation.
        \begin{align*}
            &H    _{1}^{( 1 )} =                                                                         \check{B   }_{1},  & &                    \\
            &H    _{i}^{( 1 )} =   \tilde{F}_{i}        H_{i - 1}^{( 1 )}                              + \check{B   }_{i},  & & i = 2, \dots , N,  \\
            &\Phi _{i}         =                 \left( H_{i    }^{( 1 )}                 \right) ^{2}                   ,  & & i = 1, \dots , N,  \\
            &H    _{1}^{( 2 )} =                                            \Phi _{    1}                                ,  & &                    \\
            &H    _{i}^{( 2 )} =   \tilde{F}_{i} \left( H_{i - 1}^{( 2 )} + \Phi _{i - 1} \right)      +        \Phi _{i},  & & i = 2, \dots , N.
        \end{align*}
    \end{cor}

    Now, we give an implementation.
    We compare this implementation with that in \cite{YKN12}.
    Similarly to the implementation in \cite{YKN12}, we use techniques for optimization of implementation.
    Firstly, we try to reduce the number of ``loops'' by the technique of ``loop fusion''.
    We try also to reduce the number of divisions which takes a longer time than multiplications.
    Next, we try to raise ``register hit rate'' or ``cash hit rate'' by an attempt to reduce ``working memories''.
    We avoid use of an ``array'' if it is not necessary.
    In contrast to the implementation in \cite{YKN12}, the implementation in this paper does not require an array except for the one to store the input data.
    Lastly, we try to raise ``cash hit rate'' by using the same ``variable'' consecutively.
    For comparison with the implementation in \cite{YKN12}, let $b_{i} = \sqrt{q_{i}}$ for $i = 1, \dots , N$ and $c_{i} = \sqrt{e_{i}}$ for $i = 1, \dots , N - 1$ be input data.
    Let $b_{i}$ for $i = 1, \dots , N    $ be recorded in ``array'' B$[i]$.
    Let $c_{i}$ for $i = 1, \dots , N - 1$ be recorded in ``array'' C$[i]$.
    An algorithm for computing the trace $\textrm{Tr} ( ( B B^{\top} )^{- 2} )$ based on the recurrence relation in this paper is shown in Algorithm \ref{Imp2nd_TP}.
    The variables H1, P, IB and F store the constants $H_{i}^{( 1 )}$, $\Phi _{i} = ( H_{i}^{( 1 )} )^{2}$, $\check{B}_{i}$ and $\tilde{F}_{i}$, respectively.
    The variable H2 in the line 3 stores $                H_{1    }^{( 2 )}$.
    The variable H2 in the line 8 stores $\tilde{F}_{i} ( H_{i - 1}^{( 2 )} + \Phi _{i - 1} )$ temporarily.
    After computing of $\Phi _{i}$ in the line 10, the variable H2 in the line 11 stores the constant $H_{i}^{( 2 )}$.
    The variable J is used to compute the trace.
    Our implementation requires only one loop while that in \cite{YKN12} requires two loops.
    Our implementation use no array except for the one to store the input data, in contrast to that in \cite{YKN12}.
    The numbers of arithmetic operations in our implementation and that in \cite{YKN12} are shown in Table \ref{Table_Trace2}.
    It can be seen that the former is smaller.
    Thus, it is expected that the execution time for the computation of the traces with our implementation is shorter than that with implementation in \cite{YKN12}.
    \begin{algorithm}
        \caption{An implementation of an algorithm for computing the trace $\text{Tr} ( ( B B^{\top} )^{- 2} )$ with a method based on the recurrence relation
                 in this paper} \label{Imp2nd_TP}
        \begin{algorithmic}[1]
            \STATE $\mathrm{H1} \gets 1.0 / ( \mathrm{B}[1] * \mathrm{B}[1] )$
            \STATE $\mathrm{P}  \gets         \mathrm{H1}   * \mathrm{H1}    $
            \STATE $\mathrm{H2} \gets         \mathrm{P}                     $
            \STATE $\mathrm{J}  \gets         \mathrm{H2}                    $
            \FOR{$i = 2$ \TO $N$ {\bf by} $+1$}
                \STATE $\mathrm{IB} \gets 1.0               / ( \mathrm{B}[i]     * \mathrm{B}[i] )$
                \STATE $\mathrm{F}  \gets \mathrm{C}[i - 1] *   \mathrm{C}[i - 1] * \mathrm{IB}    $
                \STATE $\mathrm{H2} \gets \mathrm{F}        * ( \mathrm{H2}       + \mathrm{P}    )$
                \STATE $\mathrm{H1} \gets \mathrm{F}        *   \mathrm{H1}       + \mathrm{IB}    $
                \STATE $\mathrm{P}  \gets \mathrm{H1}       *   \mathrm{H1}                        $
                \STATE $\mathrm{H2} \gets \mathrm{H2}       +   \mathrm{P}                         $
                \STATE $\mathrm{J}  \gets \mathrm{J}        +   \mathrm{H2}                        $
            \ENDFOR
            \RETURN $\mathrm{J}$
        \end{algorithmic}
    \end{algorithm}
    \begin{table}[h]
        \begin{center}
            \caption{Comparison of the number of arithmetic operations in computation of $\text{Tr} ( ( B B^{\top} )^{- 2} )$}  \label{Table_Trace2}
            \begin{tabular}{l | c c}
                                &  this paper  &  Ref. \cite{YKN12}  \\
                \hline
                addition        &  $4 N - 4$   &  $5 N - 5$          \\
                multiplication  &  $6 N - 4$   &  $8 N - 6$          \\
                division        &  $  N    $   &  $  N    $
            \end{tabular}
        \end{center}
    \end{table}

    \section{Concluding Remarks}  \label{ConcRmks}

    In this paper, we present a new formula for the traces of inverse powers of a positive definite symmetric tridiagonal matrix.
    From these traces, lower bounds of the minimal singular value of an upper bidiagonal matrix are obtained.
    The formula consists of only addition, multiplication and division among positive quantities, namely, it is subtraction-free.
    This property clearly excludes any possibility of cancellation error.
    Derivation of this formula is based on an idea quite different from that in \cite{YKN12}.
    
    An efficient implementation for the trace $\textrm{Tr} ( ( B B^{\top} )^{- 2} )$ which is useful in practice is also presented.
    This implementation has a few merits compared with that in \cite{YKN12}. 
    An application of the new formula to singular value computing is shown in \cite{TTIKYIN13}.

\end{document}